\documentclass[12pt]{article}
\usepackage{geometry}
\usepackage{epsfig}
\begin{document}

\def\hind{\hangindent=2pc\hangafter=1}
\newfont{\smcaps}{cmcsc10 scaled\magstep1}
\newfont{\tti}{cmss10 scaled\magstep1}
\newfont{\ttiSL}{cmssi10 scaled\magstep1}
\newfont{\ttiL}{cmss10 scaled\magstep3}
\font\tty=cmtt10 at 11truept
\newcommand{\var}{{\rm ~Var\,}}
\newcommand{\sd}{{\rm ~se\,}}
\newcommand{\sign}{{\rm ~sign\,}}
\newcommand{\Cov}{{\rm ~Cov\,}}
\newcommand{\E}{{\rm ~E\,}}
\newcommand{\BCA}{${{\rm ~BC}_a\,}$}
\newcommand{\AR}{{\rm ~AR\,}}
\newcommand{\MA}{{\rm ~MA\,}}
\newcommand{\ARMA}{{\rm ~ARMA\,}}
\newcommand{\B}{{\cal ~B\,}}

\baselineskip=22pt

\title{Fitting MA$(q)$ Models in the Closed Invertible Region}
\date{}
\author{Y. Zhang$^{a.*}$ and A.I. McLeod$^{b.1}$\\
$^{a}$Acadia University and
$^{b}$The University of Western Ontario\\
}
\maketitle
\medskip
\noindent Preprint. Zhang, Y. and McLeod, A.I. (2006).
Fitting MA(q) Models in the Closed Invertible Region,
{\it Statistics and Probablity Letters}, 76, 1331-1334. \hfill
\medskip
\hrule
\bigskip
{\bf Abstract\/}
\break
The use of reparameterization in the maximization of the likelihood function
of the $\MA(q)$ model is discussed.
A general method for testing for the presence of a parameter estimate
on the boundary of an MA$(q)$ model is presented.
This test is illustrated with a brief simulation experiment for
the MA$(q)$ for $q=1,2,3,4$ in which it
is shown that the probability of an estimate being on the
boundary increases with $q$.

\smallskip
\noindent {\it Key Words:\/}
Admissible region for the autoregressive-moving average time series;
ARMA model reparameterization;
Numerical maximum likelihood estimation.
\bigskip
\hrule
\bigskip
\noindent $^{*}$Corresponding Author: A.I. McLeod \hfill
\break\noindent Department of Statistical and Actuarial Sciences \hfill
\break\noindent The University of Western Ontario \hfill
\break\noindent London, Ontario N6A 5B7 \hfill
\break\noindent aimcleod@uwo.ca \hfill
\break\noindent $^{1}$A.I. McLeod supported by NSERC Discovery Grant.

\newpage
\noindent{\bf 1. Introduction\/}
\bigskip
\baselineskip=22pt

The $\MA(q)$ model with mean zero may be written,
$z_t = \theta(B) a_t$, where $\theta(B)=1-\theta_1 B-\ldots \theta_q B^q$
and $a_t$ is Gaussian white noise with variance $\sigma_a^2$.
The model is said to be invertible if all roots of $\theta(B)=0$
lie outside the unit circle (Box et al., 1991, \S 3.3.1).
$\MA(q)$ models, with roots inside the unit circle, may be reparameterized
so that the roots are outside the unit circle (Brockwell and Davis, 1991, p.88).
When using an exact maximum likelihood algorithm, such as the innovations
algorithm (Brockwell and Davis, 1991, \S 8.6), some parameter estimates
may lie on the non-invertible boundary (Kang, 1975).
Cryer and Ledolter (1981) proved algebraically that there is a positive
probability of the maximum likelihood estimator being on the unit
boundary in the $\MA(1)$ case.
The use of reparameterization to obtain maximum likelihood
estimates over the closed invertible region is discussed in \S 2.

A noninvertible $\MA(q)$ model means that $a_t$ cannot be
expressed in terms of past observations. For this reason, in most
situations many time series analysts would prefer an invertible
model. Additionally, with a noninvertible model care must be taken
in the algorithms used to compute the residuals and forecasts.
Statistical inference for the parameters in a noninvertible model
also becomes more difficult. For these reasons it is recommended
that a model be tested to determine if it has a parameter estimate
on the noninvertible boundary. In those cases where a parameter estimate on the
boundary is found, there are a number of remedies. One simple
approach would be to consider a different type of time series
model. For example, in \S 3  it is noted that in the mixed
$\ARMA(p,q)$ case, where $p>0$, maximum likelihood estimates on
the noninvertible boundary are less likely. Another approach would
be to use mean likelihood estimation (McLeod and Quenneville,
2000) or a Bayesian approach (Marriott and Smith, 1992). In
addition to guaranteeing invertibility, these estimation
techniques have the same first-order asymptotic efficiency as
maximum likelihood and the mean-square error of the parameters
estimates is usually less (McLeod and Quenneville, 2000). In \S 4
we present a convenient test for the presence of an estimate on
the noninvertible boundary.
It should be pointed out that one cannot simply constrain the
maximum likelihood estimates to be inside the invertible region
since from the results of Cryer and Ledolter (1981),
estimates on the boundary will cause the constrained likelihood
approach to either fail by converging to an estimate outside
the admissible region or else, if a penalty function approach
is used for the maximization, the estimates will be very close to
the boundary.
Neither of these situations is desirable.

\bigskip
\noindent{\bf 2.  Reparameterization For Constrained Maximum Likelihood Estimation\/}
\bigskip

Let ${\cal D}_\theta$ denote the invertible region for an
$\MA(q)$ in the parameter space $(\theta_1,\ldots,\theta_q)$.
As noted by Monahan (1984), the reparameterization discussed by
Barndorff-Neilsen and Schou (1973) may be extended for use with the $\ARMA(p,q)$.
For simplicity we discuss the $\MA(q)$ case but it is easy
to extend our results to the moving-average parameters in the $\ARMA(p,q)$.
Monahan (1984) defined a transformation,
${\cal B}:(\zeta_1,\ldots,\zeta_q) \longrightarrow (\theta_1,\ldots,\theta_q)$.
For brevity, let $\zeta=(\zeta_1,\ldots,\zeta_q)$,
$\theta=(\theta_1,\ldots,\theta_q)$
and
$\theta_i = {\cal B}_i(\zeta), i=1,\ldots,q$.
This transformation may be computed using the recursion,
\begin{equation}
\theta_{i,k} = \theta_{i,k-1}+\zeta_k \theta_{k-i,k-1},
\quad i=1,...,k-1; k=1,\ldots, q,
\newcounter{MonahanTransformation}
\setcounter{MonahanTransformation}{\value{equation}}
\end{equation}
where $\theta_{i,q}=\theta_i, \zeta_i=\theta_{i,i}, i=1,\ldots,q$.
The invertible region for the transformed parameters $\zeta$,
denoted by  ${\cal D}_\zeta$, is simply the interior of the unit cube,
$|\zeta_i|<1, i=1,\ldots,q$.
Barndorff-Neilsen and Schou (1973) showed that ${\cal B}$ is 1:1 and
onto as well as continuously differentiable with a continuously differentiable inverse
inside the invertible region.
However, since the determinant Jacobian of the transformation
${\cal B}$ is zero on the non-invertible boundary
(Barndorff-Neilsen and Schou, 1973, p.414), the transformation is
not 1:1 there and consequently the inverse function is not well defined.
For example in the $\MA(2)$ case,
${\cal B}(\zeta_1, \zeta_2)= (\zeta_1(1-\zeta_2),\zeta_2)$
and so  ${\cal B}(\zeta_1, 1)=(0,1)$.

Denote the likelihood function of the $\MA(q)$ by ${\cal L}(\theta)$.
Then the reparameterized likelihood function may be written
${\cal L}({\cal B}(\zeta))$,
where $\zeta$ is constrained, $|\zeta_i|\le 1, i=1,\ldots,q$.
Standard minimization algorithms, such as those implemented in
{\it Mathematica\/} (Wolfram, 2005), may be used to obtain the
maximum likelihood estimates over the closed invertible region.

The transformation ${\cal B}$ is continuous and differentiable in the closed
invertible region.
In order to use the transformation $\cal B$ for maximum likelihood
estimation in the case where the estimates may lie on the
boundary of ${\cal D}_\theta$,
it is necessary that each point on the boundary
of ${\cal D}_\theta$ be the image of a point on the boundary of ${\cal D}_\zeta$.
In other words, it is necessary that ${\cal B}$ be onto.
This is proved in \S 3.

\bigskip
\noindent{\bf 3. Proof ${\cal B}$ Is Onto\/}
\bigskip

Denote the boundary sets of ${\cal D}_\theta$ and ${\cal D}_\zeta$
by $\partial_\theta$ and $\partial_\zeta$ respectively.
Theorem 2 of Barndorff-Nielsen and Schou (1973) showed that
${\cal B}$ maps ${\cal D}_\zeta$ onto ${\cal D}_\theta$
and that this mapping is one-to-one.
But ${\cal B}$ is no longer one-to-one on the boundary $\partial_\zeta$.
It is non-trivial to show that ${\cal B}$ maps $\partial_\zeta$
onto $\partial_\theta$.
After careful investigation and discussion with an expert in point
set topology we were not able to find any suitable theorem which is
applicable to this situation
and so for completeness we have included a proof from first principles.

\noindent {\bf Theorem 1.\/}
$\partial_\theta={\cal B}(\partial_\zeta)$

\noindent{\bf Proof.\/}
First we show that ${\partial_\theta} \subset {\cal B} ({\partial_\zeta})$.
Let $\bar {\cal D}_\theta$ and $\bar {\cal D}_\zeta$
denote the closures of ${\cal D}_\theta$ and ${\cal D}_\zeta$ respectively.
Since $\cal B$ is a polynomial and hence continuous on
$\bar {\cal D}_\zeta$, ${\cal B}( \bar{\cal D}_\zeta) \subset
\overline{{\cal B}({\cal D}_\zeta)} ={\bar {\cal D}_\theta}$.
Meanwhile,
${\cal B}(\bar {\cal D}_\zeta)=
{\cal D}_\theta \cup {\cal B}(\partial_\zeta)$
is closed since $\bar {\cal D}_\zeta$ is compact and therefore
${\cal B}(\bar {\cal D}_\zeta)\supset \bar {\cal D}_\theta$.
Hence ${\cal B}(\bar {\cal D}_\zeta)=\bar {\cal D}_\theta$, so
${\cal D}_\theta \cup {\cal B}(\partial_\zeta)={\cal D}_\theta \cup \partial_\theta$.
Since ${\cal B}$ is a homeomorphism between ${\cal D}_\theta$ and
${\cal D}_\zeta$, it follows that ${\cal D}_\theta$ is open since
${\cal D}_\zeta$ is open.
${\cal D}_\theta \cap\partial_\theta = \emptyset$.
Hence ${\partial_\theta} \subset {\cal B}({\partial_\zeta})$.

Next, we show $\partial_\theta \supset {\cal
B}(\partial_\zeta)$. Let $\vartheta \in \partial_\zeta$. There
exists a sequence $\{\vartheta_n \}\in {\cal D}_\zeta$ such that
$\vartheta_n \rightarrow \vartheta$. Hence ${\cal
B}(\vartheta_n)=\vartheta_n \in {\cal D}_\theta \rightarrow {\cal
B}(\vartheta)=\vartheta$ by the continuity of ${\cal B}$ on $\bar
{\cal D}_\zeta$. If $\vartheta \in \cal D_\theta$, ${\cal B}^{-1}
({\vartheta_n})=\vartheta_n \rightarrow {\cal B}^{-1}
(\vartheta)=\vartheta \in {\cal D}_\zeta $ by the continuity of
${\cal B}^{-1}$ on ${\cal D}_\theta$, which causes a contradiction
with $\vartheta \in \partial_\zeta$. Therefore, ${\cal
B}(\vartheta)\notin \cal D_\theta$. It follows that ${\cal
B}(\vartheta) \in \partial_\theta$ since ${\cal D}_\theta \cup {\cal
B}(\partial_\zeta)={\cal D}_\theta \cup \partial_\theta$. Hence
$\partial_\theta \supset {\cal B}(\partial_\zeta)$.

\bigskip
\noindent{\bf 4. Test for Estimate on the Noninvertible Boundary\/}
\bigskip

Monahan (1984) showed that inside the invertible region,
${\cal B}^{-1}$, may be obtained using the recursive formula,
\begin{equation}
\theta_{i,k-1} =
(\theta_{i,k} - \theta_{k,k} \theta_{k-i,k})/(1-\theta_{k,k}^2), \quad i=1,\ldots, k-1;
k=q,\ldots, 1,
\newcounter{MonahanInverse}
\setcounter{MonahanInverse}{\value{equation}}
\end{equation}
where
$\theta_i = \theta_{i,q}, \zeta_i = \theta_{i,i}, i=1,\ldots,q$.
More generally, for the closed invertible region,  define
${\cal B}^-: \bar{\cal D}_\theta \longrightarrow \bar{\cal D}_\zeta$
using the recursion,
\begin{eqnarray}
\theta_{i,k-1} & = &
(\theta_{i,k} - \theta_{k,k} \theta_{k-i,k})/(1-\theta_{k,k}^2),
\quad {\rm if\/}\ |\theta_{j,j}|<1,\ {\rm for\ all\/}\ j \ge k, \nonumber \\
& = & 0 \quad {\rm otherwise\/},
\newcounter{PrincipalValue}
\setcounter{PrincipalValue}{\value{equation}}
\end{eqnarray}
where
$i=1,\ldots, k-1; k=q,\ldots, 1.$
In ${\cal D}_\theta$, ${\cal B}^-$ is the same as ${\cal B}^{-1}$
and, using mathematical induction, it may be shown that
all points on the boundary of ${\cal D}_\theta$ are mapped
into boundary points in ${\cal D}_\zeta$.
As an illustration,
$${\cal B}(\zeta _1,1,\zeta _3,\zeta _4)=
\left(-\zeta _3 \left(\zeta _4+1\right),1-\zeta _4,\zeta _3
   \left(\zeta _4+1\right),\zeta _4\right)$$
and
$${\cal B}^-(-\zeta _3 \left(\zeta _4+1\right),1-\zeta _4,\zeta _3
   \left(\zeta _4+1\right),\zeta _4)=
\left(0,1,\zeta _3,\zeta _4\right).$$

Let $\hat \theta$ denote estimates which belong to the closed invertible region
and let $\hat \zeta = {\cal B}^{-}(\hat \theta)$,
where $\hat\theta \in \bar{\cal D}_\theta$.
Then $|\hat\zeta_i|\leq 1$ and $\hat \theta$
is on the boundary if and only if
$\hat \zeta_i =\pm 1$ for at least one $i \in \{1,\dots, q\}$.
In practice, we need to take into account the finite precision and
rounding error in our computations so we may declare an estimate
is on the boundary if $|1-|\hat \zeta_i||<\epsilon$ for any
$i \in \{1,\dots, q\}$.
In the simulation example below, we took $\epsilon=10^{-6}$.
If a different software environment were used
it might be advisable to take $\epsilon$ somewhat larger.
The use of this test is illustrated in \S 5.

\bigskip
\noindent{\bf 5. Simulation Experiment\/}
\bigskip

A simulation experiment was conducted to see how the probability
of an estimate on the boundary
depends on the model order $q$ and
series length $n$. An $\MA(1)$ with parameter $\theta_1=-0.9, -0.6,
\ldots, 0.9$,  was simulated 100 times and for each
simulation the exact maximum likelihood estimates were determined
using the innovation algorithm and the
{\tty NMinimize\/} function in {\it Mathematica\/} (Wolfram, 2005)
for each of the MA$(q)$ models with $q=1,2,3,4$. An estimate was counted as on the
boundary if any one of $\hat\zeta_1, \ldots, \hat\zeta_q$ was
within $10^{-6}$ of $\pm 1$. The results are summarized in Table 1.
We see that the probability of an estimate on the boundary
increases with the model order $q$, decreases with sample size $n$
and increases with the distance of the true parameter to the unit
boundary. Note that the test for an estimate on the boundary
requires Theorem 1.
For the model MA$(1)$, the results agree well
with the previous studies mentioned above.

In another simulation experiment, $\ARMA(1,q), q=1,2,3,4$
models were fit to a series of length 25
generated by an $\ARMA(1,1)$ model with parameters
$\phi=-0.5, 0, 0.5$ and $\theta_1=-0.9, -0.6, \ldots, 0.9$.
It was found that the probability of a root on the moving average boundary
was quite small but appeared to be nonzero in all cases.

\bigskip
\noindent{\bf 6. Concluding Remarks\/}
\bigskip

Theorem 1 is also useful in the maximum likelihood estimation of
$\ARMA(p,q)$ models, $\phi(B) z_t = \theta(B) a_t$ where
$\phi(B)=1-\phi_1 B-\ldots \phi_p B^p$
and
$\theta(B)=1-\theta_1 B-\ldots \theta_q B^q$.
For a causal-stationary process, $\phi(B)=0$ has all roots outside the
unit circle and for model identifiability, in the econometric sense (Harvey, 1990, \S 3.6),
we require that $\theta(B)=0$ has all roots on or inside the unit circle.
The results \S 2 and \S 4 may be extended to the case.
Complete details of the maximum likelihood algorithm for $\ARMA(p,q)$ models,
as well as of the simulations
in \S 5, are given in {\it Mathematica\/} notebooks available from the authors.

\newpage
\noindent{\bf Acknowledgements\/}
\bigskip

The authors would like to thank Dr. Andr\'e Boivin for
his suggestions and comments on Theorem 1.

\bigskip
\noindent{\bf References\/}
\bigskip

\noindent Barndorff-Nielsen, O. and Schou G. (1973),
On the Parametrization of Autoregressive Models by Partial Autocorrelations,
{\it Journal of Multivariate Analysis\/} {\bf 3}, 408--419.

\noindent  Box, G.E.P., Jenkins, G.M. Reinsel, G.C., (1994),
{\it Time Series Analysis: Forecasting and Control\/}
(Prentice-Hall, Englewood Cliffs, 3rd ed.).

\noindent Brockwell, P.J. and Davis, R.A. (1991),
{\it Time Series: Theory and Methods\/}
(Springer-Verlag, New York, 2nd ed.).

\noindent Cryer, J. D. and Ledolter, J. (1981),
Small-sample properties of the maximum likelihood estimator in the first-order
moving average model,
{\it Biometrika\/} {\bf 68}, 691--694.

\noindent Hamilton, J.D. (1994),
{\it Time Series Analysis\/}
(Princeton University Press, Princeton).

\noindent Harvey, A. (1990),
{\it The Econometric Analysis of Time Series\/}
(MIT Press, Cambridge, 2nd ed.).

\noindent Kang, K.M. (1975). A Comparison of Estimators for Moving
Average Proesses. Unpublished manuscript, Australian Bureau of
Statistics.


\noindent Marriott, J. M. and Smith, A. F. M. (1992),
Reparametrization Aspects of Numerical Bayesian Methodology for
Autoregressive Moving-average Models,
{\it Journal of Time Series Analysis\/} {\bf 13}, 327--343.

\noindent McLeod, A.I. and Quenneville, B. (2000),
Mean likelihood estimators,
{\it Statistics and Computing\/} {\bf 11}, 57--65.

\noindent Monahan, J.F. (1984),
 A note on enforcing stationarity in autoregressive-moving average models,
{\it Biometrika\/} {\bf 71}, 403--404.

\noindent Wolfram, S. (2005),
{\it The Mathematica Book 5\/}
(Wolfram Media, Champaign).

\strut
\begin{table}
\caption{Proportion of times an estimate was on the noninvertible boundary
in 100 simulations of a series of length $n$ from an $\MA(1)$ with parameter
$\theta_1$ when it is fit with an $\MA(q), \quad q=1,2,3,4$. }
$$\vbox{\halign{\strut
 \hfill#\quad    
&\hfill#\quad    
&\hfill#\quad    
&\hfill#\quad    
&\hfill#\quad    
&\hfill#\quad    
&\hfill#\quad    
&\hfill#\quad    
&\hfill#\quad    
&\hfill#\quad    
&\hfill#\quad    
&\hfill#\quad    
&\hfill#\quad    
&\hfill#\quad    
\cr             
\noalign{\hrule}
\noalign{\smallskip}
\noalign{\hrule}
\noalign{\smallskip}
$n$&$\theta$&$\MA(1)$&$\MA(2)$&$\MA(3)$&$\MA(4)$
\cr
\noalign{\smallskip}
\noalign{\hrule}
\noalign{\smallskip}
$25$&$-0.9$&$0.53$&$0.51$&$0.51$&$0.46$
\cr
$25$&$-0.6$&$0.08$&$0.15$&$0.23$&$0.53$
\cr
$25$&$-0.3$&$0.00$&$0.03$&$0.17$&$0.27$
\cr
$25$&$0.0$&$0.01$&$0.05$&$0.24$&$0.35$
\cr
$25$&$0.3$&$0.02$&$0.06$&$0.23$&$0.27$
\cr
$25$&$0.6$&$0.14$&$0.21$&$0.37$&$0.39$
\cr
$25$&$0.9$&$0.49$&$0.51$&$0.49$&$0.57$
\cr
$50$&$-0.9$&$0.36$&$0.36$&$0.37$&$0.32$
\cr
$50$&$-0.6$&$0.01$&$0.03$&$0.06$&$0.11$
\cr
$50$&$-0.3$&$0.00$&$0.03$&$0.04$&$0.08$
\cr
$50$&$0.0$&$0.00$&$0.00$&$0.02$&$0.04$
\cr
$50$&$0.3$&$0.01$&$0.02$&$0.02$&$0.04$
\cr
$50$&$0.6$&$0.01$&$0.03$&$0.03$&$0.07$
\cr
$50$&$0.9$&$0.26$&$0.33$&$0.28$&$0.36$
\cr
\noalign{\smallskip}
\noalign{\hrule}
}}$$
\end{table}

\end{document}